\documentclass[12pt]{article}

\usepackage{times}
\usepackage[centertags]{amsmath}
\usepackage{amssymb}
\usepackage{amsthm}
\usepackage[latin1]{inputenc}
\usepackage[colorlinks=true,urlcolor=red,citecolor=red,linkcolor=red]{hyperref}
\usepackage[english]{babel}
\usepackage{graphicx}
\usepackage{tikz}
\usepackage{qtree}
\usepackage{tikz-qtree}
\usepackage[lined,algonl,boxed]{algorithm2e}
\usetikzlibrary{trees,positioning,shapes} 
\newtheorem{prop}{Proposition}[section]
\newtheorem{raff}[prop]{Proposition}

\newtheorem{Lemma3}[prop]{Lemma}

\newtheorem{genere}[prop]{Theorem}
\theoremstyle{definition} 
\newtheorem{example1}[prop]{Example}
\newtheorem{example2}[prop]{Example}
\newtheorem{example3}[prop]{Example}
\newtheorem{example4}[prop]{Example}
\newtheorem{example44}[prop]{Example}
\newtheorem{example5}[prop]{Example}
\newtheorem{example6}[prop]{Example} 

\newtheorem{example8}[prop]{Example}
\newtheorem{remark0}[prop]{Remark}
\newtheorem{remark1}[prop]{Remark}
\newtheorem{remark2}[prop]{Remark}

\newtheorem{remark4}[prop]{Remark}

    {\cleardoublepage\thispagestyle{empty}\null\vfill\begin{center}%
    \bfseries Acknowledgements\end{center}}%
    {\vfill\null}

\newcounter{lastnote}

\title{Arf good semigroups with fixed genus}

\author
{Giuseppe Zito\\
\\
\\
}

\date{}

\begin{document}

\baselineskip16pt

\maketitle

\begin{abstract}
In this paper we give an algorithm for the computation of all the Arf numerical semigroups with a given genus. Moreover, we  generalize the concept of genus of a numerical semigroup to good semigroups of $\mathbb{N}^r$ and we give a procedure to calculate all the Arf semigroups of $\mathbb{N}^r$ with a given genus.
\end{abstract}

\let\thefootnote\relax

\section*{Introduction}

In this paper we  present  procedures concerning  Arf subsemigroups of $\mathbb{N}^r$ for $r \geq 1$.  In particular, for $r=1$, we deal with the concept of Arf numerical semigroup whose importance is motivated by its role in the  study of equivalence between  algebroid branches. Specifically, two algebroid branches are said to be equivalent if they have the same multiplicity sequence (cf. \cite[Definition 1.5.11]{Campi}), that is the sequence of the multiplicities of their succesive blow ups. In \cite{Arf}, Arf gave an algebraic interpretation of this equivalence,  introducing the concept of Arf ring and showing that for each ring $R$ there exists a smallest Arf overring $R'$, called the Arf closure of $R$. A ring $R$ and its Arf closure have the same multiplicity sequence, and from this it follows that two algebroid branches are equivalent if and only if their Arf closure have the same value semigroup. Taking in account the properties of the value semigroup of an Arf ring,  a numerical Arf semigroup is defined as a numerical semigroup $S \subseteq \mathbb{N}$ such that $S(s)-s$ is a semigroup, for each $s\in S$, where $S(s)=\left\{n \in S; n\geq s \right\}$.
The multiplicity sequence of an Arf numerical semigroup $S=\left\{ s_0=0<s_1<s_2,\ldots \right\}$ is very easy to compute since it coincides with the sequence of the differences $s_{i+1}-s_i$.  By focusing on the numerical properties that a sequence has to satisfy to be a multiplicity sequence, it is possible to study the Arf numerical semigroups  with a combinatorial approach without referring to the ring theory
context.
 In the first section of this paper we present a new procedure to compute all Arf numerical semigroups $S$  with a prescribed  genus, that is the cardinality of $\mathbb{N} \setminus S$ (this problem was already addressed in \cite{algo}). Then we start dealing with the case $r>1$, considering the class of Arf good semigroups of  $\mathbb{N}^r$. The definition of good semigroup was presented  \cite{BDF} and it is based on the properties of the value semigroups of one dimensional analytically unramified ring (for example the local ring of an algebraic curve), although
  it is possible to show that the resulting class of good semigroups is  bigger than the class of value semigroups. Therefore the good semigroups can be seen as a natural generalization of the numerical semigroups.   A good semigroup $S \subseteq \mathbb{N}^n$ is said to be local if the zero vector is the only element of $S$ with zero coordinates.

Then we can give  the definition  of Arf semigroup in the context of  good semigroups of $\mathbb{N}^r$, by adapting in a natural way the numerical version, considering the usual partial ordering of $\mathbb{N}^r$.  Furthermore, in order to extend  the concept of multiplicity sequence, in \cite{BDF} it is also shown that to each local Arf good semigroup, it can be associated a multiplicity tree that characterizes the semigroup completely. A tree $T$ of vectors of $\mathbb{N}^n$ has to satisfy some properties to be a multiplicity tree of a local Arf good semigroup. For instance, it must have multiplicity sequences along its branches (because its projections are Arf numerical semigroups) and each node must be able to be expressed as a sum of nodes in a subtree of $T$ rooted in it.
In \cite{Zito}, taking in account this $1$-$1$ correspondence, it is showed a  way to describe a multiplicity tree of an Arf semigroup of  $\mathbb{N}^r$  by mean of an ordered collection $E$ of multiplicity sequence and a matrix $M(T)_E=(p_{i,j})$. In the same paper it is shown that we can consider a particular class of trees, the untwisted ones, that are easier to study because they can be expressed by a vector of $\mathbb{N}^{r-1}$ instead of a matrix. This is not restrictive because every multiplicity tree can be transformed by an opportune permutation into an untwisted one (this implies that the corresponding semigroups are isomorphic).

  The conductor $\textbf{c}$ of a good semigroup is the minimal vector  such that $ \textbf{c}+ \mathbb{N}^r \subseteq S$ (the existence of such a vector is guaranteed by the properties of good semigroups). In \cite{Zito2},  it is shown how to deduce the value of the conductor of an Arf semigroup from the representation of its multiplicity tree and the author presented procedures for the computation of all Arf semigroups of $\mathbb{N}^r$ with a prescribed conductor.
In order to generalize the procedure introduced in the numerical case,  we need to extend the concept of genus  to good semigroups of $\mathbb{N}^r$. This can be done by using the function $d(\textunderscore\setminus \textunderscore )$ defined in \cite{D'a} .
In particular, if $S$ is a good semigroup of $\mathbb{N}^r$ with conductor $\delta$, then $C=\delta+\mathbb{N}^r$ is an ideal  of $S$ and we can define the genus of $S$ as
$ g(S)=d(\mathbb{N}^r \setminus C) - d(S \setminus C). $
 Then, using the aforementioned representation for multiplicity trees and the procedure developed in the numerical case, we address the problem of computing the set $\textrm{Gen}(r,n)$ of all the Arf good semigroups of $\mathbb{N}^r$ with an untwisted tree and genus $n$.

The structure of the paper is the following.

In Section \ref{section4}, we firstly recall some definitions and properties concerning the numerical semigroups. Then we address the problem, already studied in \cite{algo}, of finding  the set $\textrm{Gen}(n)$ of the multiplicity sequences of all the Arf numerical semigroups with  genus $n$. We accordingly adapt the procedure presented in \cite{Zito2} for the computation of all the Arf numerical semigroups with a given conductor,  presenting a non recursive algorithm for the computation of $\textrm{Gen}(n)$, that appears to be faster than the one presented in \cite{algo}.

In Section \ref{section5}, we begin to deal with semigroups in higher dimensions. In order to do that we recall the definition of good semigroups and the properties of the multiplcity tree of an Arf local good semigroup and we accordingly generalize the concept of genus to good semigroups of $\mathbb{N}^r$.
  We define the sets $\textrm{Gen}(r,n)$ of the untwisted multiplicity trees of all the Arf good semigroups of $\mathbb{N}^r$ with genus $n$. In Theorem \ref{genere} we give a way to compute the genus of an untwisted Arf good semigroup of $\mathbb{N}^r$ from its representation $T_E$. Using this theorem it is possible to develope, using induction on $r$, a procedure that calculate the set $\textrm{Gen}(r,n)$, and from Proposition \ref{raff} we deduce a possible way to make this computation faster. We give a strategy for computing the set $\overline{\textrm{Gen}(r,n)}$ of all the possible multiplicity trees (twisted and untwisted) associated to an Arf semigroup with genus $n$ in $\mathbb{N}^r$. At the end of the section we give an example of the application of the developed procedure and we present some tables containing  cardinalites of  the sets $\textrm{Gen}(r,n)$ for some values of $r$ and $n$. 

The procedures presented here have been implemented in GAP (\cite{gap}).
\section{An algorithm for the Arf numerical semigroups with a given genus} \label{section4}
In \cite{algo} it is presented an algorithm for the computation of the set of the Arf numerical semigroups with a given genus.
In this section we give a new procedure for the computation of such a set that appeared to be faster when implemented in GAP.

First of all we recall some definitions and we fix some notations. \\ 
A numerical semigroup $S$ is a submonoid of $(\mathbb{N},+)$ such that $ \mathbb{N} \setminus S$ is finite. The maximum $F(S)$ of the set $ \mathbb{Z} \setminus S$ is known as the Frobenius number of $S$. The conductor $c(S)$ of $S$ is the smallest number such that $n \in S$ for all $n \geq c(S)$, and it is clear that we have $c(S)=F(S)+1$.

We are interested in a particular class of numerical semigroups, which is the set of Arf numerical semigroups. A numerical semigroup $S$ is said to be Arf if $S(s)-s$ is a semigroup for each $s \in S$, where $S(s)=\left\{ n \in S; n \geq s \right\}$. A nonincreasing sequence $M=\left\{m_n; n\geq 1 \right\}$ of $\mathbb{N}^*=\mathbb{N} \setminus \left\{ 0\right\}$  is called a multiplicity sequence if
\\
$\bullet\ $ There exists $k \in \mathbb{N}^*$ such that $m_n=1$ for all $n \geq k$;
\\
$\bullet\ $ For all $n \in \mathbb{N}^*$ there exists $l(n) \geq n+1$ such that $ \displaystyle m_n=\sum_{k=n+1}^{l(n)}{m_k}$.

In the following we describe a multiplicity sequence $M$ by the finite vector $ M=\nolinebreak[1] 
[m_1,\ldots,m_k]$ where $m_k$ is the last entry in $M$ different from one. We will make an exception for the costant multiplicity sequence $M=\left\{1,\ldots,1,\ldots\right\}$ that will be represented by the vector $M=[1]$. 
Given a vector $\textrm{v}$ we  always denote by $\textrm{v}[i]$ its $i$-th component

In \cite[Corollary 39]{Rosales} it is proved that a nonempty subset of $\mathbb{N}$ is an Arf numerical semigroup if and only if there exists a multiplicity sequence $M=[m_1,\ldots,m_k]$ such that $$S=\left\{0,m_1,m_1+m_2,\ldots,m_1+\dots+m_k,\rightarrow \right\},$$
where with $n \rightarrow$ we  mean that all integers larger than $n$ are in $S$.
Thus we have a $1$-$1$ correspondence between Arf numerical semigroups and multiplicity sequences and in the following we  deal mainly with the latter.
If $M=[m_1,\ldots,m_k]$ is a multiplicity sequence, we will call  $k$ the length of  $M$ and we will denote it  by $l(M)$.
We also denote by $\textrm{AS}(M)$ the Arf numerical semigroup associated to $M$. Now, from the previous result, it is clear that if $M$ is a mutliplicity sequence different from $[1]$, then the conductor of  $\textrm{AS}(M)$  is
$ \sum_{i=1}^{l(M)}{m_i},$
while the conductor of $\textrm{AS}( [1])=\mathbb{N}$ is $0$.
The genus of a numerical semigroup $S$ is the cardinality of $\mathbb{N} \setminus S$.

We denote by $\textrm{Gen}(n)$ the set the multiplicity sequences of  the Arf numerical semigroups with genus $n$.
If $M$ is a multiplicity sequence we denote by $c(M)$ the conductor of  the Arf semigroup $\textrm{AS}(M)$ associated to $M$.
It is easy to show that, if $M$ is a multiplicity sequence distinct from $[1]$, then the genus of $\textrm{AS}(M)$   is $c(M)-l(M)$. If $M=[1]$ then the genus of $\textrm{AS}(M)=\mathbb{N}$ is zero.

Suppose that we want to compute $\textrm{Gen}(n)$  with $ n \in \mathbb{N} $.
 If $n=0$ then  $\textrm{Gen}(n)=\left\{[1] \right\}$. Thus we suppose $n \geq 1$.
Denote by $$U^n(i)=\left\{ M \in \textrm{Gen}(i): M[1]+i-1 \leq n \right\}  \textrm{ for all } i=1,\ldots,n-1.$$
Now suppose that $M=[m_1,\ldots,m_k] \in \textrm{Gen}(n)$. If $k=1$ then $M=[n+1]$, otherwise we have the following situation:
\begin{itemize}
\item $2 \leq m_1 \leq n$; 
\item $c([m_2,\ldots,m_k])=c([m_1,m_2,\ldots,m_k])-m_1$ and $l([m_2,\ldots,m_k])=l([m_1,\ldots,m_k])-1$.

 So $c([m_2,\ldots,m_k])- l([m_2,\ldots,m_k])=n-m_1+1$, and
$[m_2,\ldots,m_k] \in \textrm{Gen}(n-m_1+1); $ 
\item $ m_1 \in \textrm{AS}([m_2,\ldots,m_k]);$
 \item $ m_2-m_1 \leq 0 \Rightarrow m_2+(n-m_1+1) -1 \leq n \Rightarrow [m_2,\ldots,m_k] \in  U^n(n-m_1+1)$.
\end{itemize}
So if we know the $U^n(i)$  for $i=1,\ldots,n-1$, then we can compute $\textrm{Gen}(n)$ in the following way:
$$ \textrm{Gen}(n)=\bigcup_{i=1}^{n-1}{\left\{[n-i+1, \textbf{j}]: M \in U^n(i), n-i+1 \in \textrm{AS}(\textbf{j})\right\}} \cup \left\{  [n+1] \right\}.  $$
Thus we need  a way to compute $U^n(i)$.
Suppose that $M=[m_1,\ldots,m_k] \in U^n(i)$. If $k=1$, and $i+1+i-1 =2i\leq n$ then $M=[i+1]$, otherwise we have the following situation:
\begin{itemize}
\item $2 \leq m_1 \leq i$, and $[m_2,\ldots,m_k] \in \textrm{Gen}(i-m_1+1)=\textrm{Gen}(q); $
\item $m_2+q-1=m_2+i-m_1+1-1=m_2-m_1+i \leq i \leq n \Rightarrow [m_2,\ldots,m_k] \in U^n(q); $
\item $ m_1 \in \textrm{AS}([m_2,\ldots,m_k]);$
\item $m_1+i-1 \leq n \Rightarrow 2m_1 \leq n-i+1 +m_1=n-(i+1-m_1)+2=n-q+2 \Rightarrow  m_1 \leq \nolinebreak \left \lfloor{\frac{n-q+2}{2}}\right \rfloor.$
\end{itemize}
So each $ U^n(i)$ can be constructed using $U^n(q)$ with $q<i$.
Thus we have the following algorithm for the computaton of $\textrm{Gen}(n)$.
\begin{algorithm}
\SetKwData{Left}{left}
\SetKwData{This}{this}
\SetKwData{Up}{up}
\SetKwFunction{Union}{Union}
\SetKwFunction{FindCompress}{FindCompress}
\SetKwInOut{Input}{input}
\SetKwInOut{Output}{output}
\caption{}
\Input{An integer $n$}
\Output{The set $\textrm{Gen}(n)$ of all the multiplicity sequences of Arf semigroups with genus $n$}
\BlankLine
$\textrm{Gen}(n) \longleftarrow \left\{ [n+1] \right\}$

\For{$i\leftarrow 1$ \KwTo $n-1$} {\If{$i\leq  \left \lfloor{\frac{n}{2}}\right \rfloor$}{$U^n(i) \longleftarrow \left\{ [i+1] \right\}$}\Else{$U^n(i) \longleftarrow \emptyset$}}
\For{$i\leftarrow 1$ \KwTo $n-1$}{
\For{$ M \in U^n(i)$}{
\If{$n-i+1\in \textrm{AS}(M)$}{
$\textrm{Gen}(n) \longleftarrow \textrm{Gen}(n) \cup \left\{ [n-i+1,M] \right\}$
}
\For{$k \in \textrm{AS}(M) \cap \left\{ 2,\ldots,\left \lfloor{\frac{n-i+2}{2}}\right \rfloor \right\}$} {$U^n(i+k-1) \longleftarrow U^n(i+k-1) \cup \left\{ [k,M]\right\}$}
}
}
$\textrm{Gen}(n)$
\label{algo_disjdecomp}
\end{algorithm}
\newpage
\section{Arf semigroups of $\mathbb{N}^r$ with given genus} \label{section5}
From this section we begin to deal with Arf good semigroups of $\mathbb{N}^r$.  First of all we recall some important definitions and results.
A good semigroup $S$ of $\mathbb{N}^r$ is a submonoid of $\left( \mathbb{N}^r,+\right)$ such that  (cf. \cite{BDF})\begin{itemize}
	\item For all $a,b \in S$, $	\min(a,b) \in S$; 
	\item If $a,b\in S$ and $a[i]=b[i]$ for some $i \in \left\{ 1, \ldots,n\right\}$, then there exists $c \in S$ such that $c[i] > a[i]=b[i]$, $c[j]\geq \min(a[j],b[j])$ for $j \in  \left\{ 1, \ldots,n\right\} \setminus \left\{ i\right\}$ and $c[j]=\min(a[j],b[j])$ if $a[j]\neq b[j]$;
	\item There exists $\delta \in S$ such that $\delta+\mathbb{N}^r \subseteq S$
\end{itemize}
(where we are considering the usual partial ordering in $\mathbb{N}^r$: $a \leq b$ if $a[i] \leq b[i]$ for each $i=1,\ldots,r$). The previous properties guarantee that for a good semigroup $S$ there exists a minimal vector $c(S)$ such that $c(S)+ \mathbb{N}^r \subseteq S$. Such a vector is said to be the conductor of $S$.
\\
In this paper we will always deal with local good semigroups. A good semigroup $S$  is local if the zero vector is the only vector of $S$ with some component equal to zero.
However, it can be shown that every good semigroup is the direct product of local semigroups (cf. \cite[Theorem 2.5]{BDF}).

An Arf semigroup of $\mathbb{N}^r$, is a good semigroup such that $ S(\alpha)-\alpha $ is a semigroup, for each $\alpha \in S$ , where $ S(\alpha)=\left\{\beta \in S; \beta \geq \alpha \right\}$.
The  multiplicity tree $T$ of a local Arf semigroup $S\subseteq \mathbb{N}^r$ is a tree where the nodes are vector $ \textbf{n}_i^j \in \mathbb{N}^r$, where with $\textbf{n}_i^j$ we mean that  this node is in the $i$-th branch on the $j$-th level (the root of the tree is $\textbf{n}_{1}^1=\textbf{n}_i^1 $ for all $i$ because we are in the local case and at level one all the branches must be glued) and we have 

$$ S=\left\{\textbf{0}\right\} \bigcup_{T'} \left\{ \sum_{\textbf{n}_i^j \in T' } {\textbf{n}_i^j}\right\},$$

where $T'$ ranges over all finite subtree of $T$ rooted in   $\textbf{n}_1^1$.
 
 \begin{example5}

 	Let us consider the following subset of $\mathbb{N}^2$,
 	$$ S=\left\{ (0, 0), (4, 4), (8, 6),  (12, 6)  \right\} \cup \left\{(8,8+n),(12,8+n),(14+n,6); n \in \mathbb{N} \right\} \cup $$ 
 	$$ \cup \left\{(14+m,8+n); m,n\in \mathbb{N} \right\}.$$
 	
 	\newpage
 	It is possible to verify that $S$ ia an Arf good semigroup with the following multiplicity tree:
 	
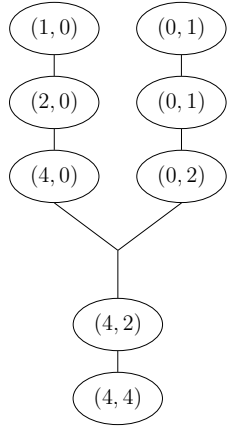
\begin{center}
 	\begin{tikzpicture}[grow'=up,sibling distance=20pt,scale=.70]
 	\tikzset{level distance=40pt,every tree node/.style={draw,ellipse}} \Tree [ .$(4,4)$ [ .$(4,2)$ [  [ .$(4,0)$  [ .$(2,0)$  [ .$(1,0)$  ] ]  ] [ .$(0,2)$  [ .$(0,1)$   [ .$(0,1)$  ] ] ]  ] ] ] ]; \end{tikzpicture} 
\end{center}
 \end{example5}

A tree $T$ is a multiplicity tree of an Arf semigroup if and only if its  nodes  satisfy the following properties (cf. \cite[Theorem 5.11]{BDF} ).

\begin{itemize}
	\item There exists $L \in \mathbb{N}$ such that for $m \geq L$, $\textbf{n}_i^m=(0,\ldots,0,1,0\ldots,0)$ (the nonzero coordinate is in the $i$-th position) for any $i=1,\ldots,n$;
	\item $\textbf{n}_i^j[h]=0$ if and only if $\textbf{n}_i^j$ is not in the $h$-th branch of the tree;
	\item Each $\textbf{n}_i^j$ can be obtained as a sum of nodes in a finite subtree $T'$ of $T$ rooted in $\textbf{n}_i^j$.
\end{itemize}
Notice that from these properties it follows that we must have multiplicity sequences along each branch. Therefore a multiplicity tree $T$ of an Arf semigroup of $\mathbb{N}^r$ can be represented by an ordered collection of $r$ multiplicity sequences $E$ and by an upper triangular matrix $ r \times r$  $$M(T)_{E}=\left(  \begin{matrix} 0 & p_{1,2} & p_{1,3} & \ldots & p_{1,r} \\ 
0 & 0 & p_{2,3} & \ldots & p_{2,r} \\ \ldots & \ldots & \ldots & \ldots & \ldots \\ 0 & 0 & 0 & \ldots& p_{r-1,r} \\ 0 & 0 & 0 & \ldots & 0 \\  \end{matrix}\right)$$where $p_{i,j}$ is the highest level such that the  $i$-th and the $j$-th branches are glued in $T$. Suppose that $E=\left\{M_1,\ldots,M_r \right\}$, then in \cite[Proposition 1.2]{Zito}, it is shown that the values that can be assigned to the $p_{i,j}$, in order to have a tree compatible with the previous properties, depend only on $M_i$ and $M_{j}$. In particular, if $M_i=M_{j},$ $p_{i,j}$ can be any positive integer,  while, if $M_i \neq M_{j}$, there exists a level $k_E(i,j) \in \mathbb{N}^*$ such that $p_{i,j}$ is at most $k_E(i,j)$.
If we set $s_{i,k}$ as the integer such that  $$ M_i[k]=\sum_{l=k+1}^{s_{i,k}}{M_i[l]},$$
then we have
$$k_E(i,j)=\min\left\{ \min( s_{i,k}, s_{j,k}): s_{i,k}\neq s_{j,k} \right\}. $$ 
Thus, if $M_1$ and $M_2$  are two distinct multiplicity sequences, we  define the compatibility between $M_1$ and $M_2$ as the integer  $\textrm{Comp}(M_1,M_2)=k_E(1,2) $ where $E=\left\{M_1,M_2 \right\}$. While if $M_1= \nolinebreak M_2$ we set by definition  $\textrm{Comp}(M_1,M_2)=+\infty$. 
A tree $T$ is untwisted if two non-consecutive branches are glued at level $l$ if and only if  all the consecutive branches  between them are glued  at a level greater or equal to $l$. We will call twisted a tree that it is not untwisted.
From the definition it follows that the matrix of an untwisted tree $T \in \tau(E)$ is such that: $$ p_{i,j}=\min\left\{ p_{i,i+1},\ldots,p_{j-1,j}\right\} \textrm{ for all } i<j. $$
So an untwisted tree can be completely described by the second diagonal of its matrix. Thus, in the following we   indicate an untwisted tree by a vector $T_E=(p_1,\ldots,p_{r-1})$ where $p_i=p_{i,i+1}$.
\begin{remark0} \label{remark0}
	It is easy to see that a twisted tree can be converted to an untwisted one by accordingly permuting its branches (the corresponding Arf semigroups are therefore isomorphic). 
	Thus in the following we can focus only on the properties of the untwisted trees, that are easier to study than the twisted ones. 
\end{remark0}
The aim of this section is to find a way to determine all the Arf good semigroups of $\mathbb{N}^r$ with a given genus.
For the Remark \ref{remark0}, we can focus only on the untwisted trees.

First of all we recall how we can extend the concept of genus for good semigroups of $\mathbb{N}^r$.
If $S$ is a numerical semigroup with conductor $c$, then $C=\left\{c,\rightarrow \right\}$ is an ideal of $S$. In this case the genus of $S$ is:
$$ g(S)=|\mathbb{N} \setminus C|-|S \setminus C|=c-|S \setminus C|.$$
So we have a natural way to extend this concept to the good semigroups of $\mathbb{N}^r$, using the function $d(\textunderscore\setminus \textunderscore )$ defined in \cite{D'a} .
If $S$ is a good semigroup of $\mathbb{N}^r$ with conductor $\delta$, then $C=\delta+\mathbb{N}^r$ is an ideal  of $S$ and we can define the genus of $S$ in the following way:
$$ g(S)=d(\mathbb{N}^r \setminus C) - d(S \setminus C). $$
Since $d(\mathbb{N}^r \setminus C)$ is the length of a saturated chain in $\mathbb{N}^r$ from the vector $ \textbf{0} \in \mathbb{N}^r$ to the conductor $\delta=(c[1],\ldots,c[r])$ (a chain $\alpha=\alpha_0<\alpha_1<\dots<\alpha_n=\beta$ in $S$ is said to be saturated if it cannot be extended to a longer chain between $\alpha$ and $\beta$ in $S$), it is easy to show that $$d(\mathbb{N}^r \setminus C)=\sum_{k=1}^r{c[k]}.$$
On the other hand, $d(S \setminus C)$ is the length of a sautared chain in $S$ from   $ \textbf{0} \in S$ to $\delta \in S$.
\begin{remark4} \label{remark4}
	Because the conductor $\delta$ is an element of the Arf semigroup $S(T)$, it can be expressed as a sum of  nodes in a subtree $T'$ of $T$. It is easy to verify that $T'$ must be the subtree consisting of  the nodes of $T$ that are different from the canonical vectors $e_i=(0,\ldots,0,1,0,\ldots,0).$
\end{remark4}

For the Arf good local semigroups with untwisted multiplicity tree we have the following theorem:

\begin{genere} \label{genere}
Suppose that $T=T_E=(p_1,\ldots,p_{r-1})$ is an untwisted mulitplicity tree of an Arf semigroup where $E=\left\{ M_1,\ldots,M_r \right\}$ is a collection of multiplicity sequences. 

Then
$$ g(S(T))=\sum_{k=1}^r{g(\textrm{AS}(M_k))}+\sum_{k=1}^{r-1}{p_k}, $$
where $S(T)$ is the Arf semigroup associated to the tree $T$ and  $\textrm{AS}(M_k)$ is the Arf numerical semigroup associated to the multiplicity sequence $M_k$.
\end{genere}
\noindent \textbf{Proof}. Denote, as usual, by $\delta=(c[1],\ldots,c[r])$ the conductor of $S(T)$, and by $C=\delta+\mathbb{N}^r$. We know that
$$ g(S(T))=d(\mathbb{N}^r \setminus C) - d(S \setminus C). $$
We have  
$$ d(\mathbb{N}^r \setminus C) =\sum_{k=1}^r{c[k]}=\sum_{k=1}^{\max(l(M_1),p_1)}{M_{1}[k]}+\ldots+\sum_{k=1}^{\max(l(M_i),p_i,p_{i-1})}{M_{i}[k]}+\ldots+\sum_{k=1}^{\max(l(M_r),p_{r-1})}{M_{r}[k]}= $$
$$ =\sum_{k=1}^{l(M_1)}{M_{1}[k]}+\max(l(M_1),p_1)-l(M_1)+\ldots+\sum_{k=1}^{l(M_i)}{M_{i}[k]}+\max(l(M_i),p_i,p_{i-1})-l(M_i)+\ldots $$ $$\ldots+\sum_{k=1}^{l(M_r)}{M_{r}[k]}+ \max(l(M_r),p_{r-1})-l(M_r)=\sum_{j=1}^{r}\left(\sum_{k=1}^{l(M_j)}{M_{j}[k]}-l(M_j)\right)+ \max(l(M_1),p_1)+$$ $$+ \sum_{k=2}^{r-1}{\max(l(M_k),p_k,p_{k-1})}+\max(l(M_r),p_{r-1}),$$
where we are using the fact that
$$ c[1]=\sum_{k=1}^{\max(l(M_1),p_1)}{M_{1}[k]}, \ldots, c[i]=\sum_{k=1}^{\max(l(M_i),p_i,p_{i-1})}{M_{i}[k]}, \ldots, c[r]= \sum_{k=1}^{\max(l(M_r),p_{r-1})}{M_{r}[k]}, $$
that was proved in \cite[Proposition 2.2]{Zito2}.

Now we want to compute $d(S \setminus C)$. We need a saturated chain in $S(T)$ from $\textbf{0}$ to $\delta$.
Suppose that we have $$ \textbf{0}=s_0<s_1<\dots<s_l=\delta,$$
a saturated chain in $S(T)$.
We clearly have $s_1=(M_{1}[1],\ldots,M_{r}[1])$, that is the multiplicity vector of $S(T)$.
Let us consider $s_q \in S(T)$, with $q=1,\ldots,l-1$. From the properties of the multiplicity tree of an Arf semigroup, there exists a subtree $T'$ of $T$, rooted in the node corresponding to $s_1$, such that $s_q$ is the sum of all the nodes belonging to $T'$.
As usual we denote by $\textbf{n}_{i}^{j}$ the node of $T$ that is in the $i$-th branch and on the $j$-th level. We denote, given a subtree $T'$ of $T$, by $\textrm{N}(T')$ the set of nodes that appears in $T'$.

Now, it is clear that, in order to have a saturated chain, $s_{q+1}$ must be the sum of all the nodes belonging to a subtree $T''$ of $T$ such that:
\begin{itemize}
\item $ T' \subseteq T''$;
\item $N(T'') \setminus N(T')=\left\{ \textbf{n}_{i}^{j} \right\}$;
\item $\textbf{n}_{i}^{j-1}\in N(T')$;
\item $\textbf{n}_{i}^{j}\neq e_i$, where $e_i $ is the $i$-th canonical vector of $\mathbb{N}^r$ (by Remark \ref{remark4} since $s_q \neq \delta$).
\end{itemize}
From the previous remark, it easily follows that
$$ d(S \setminus C)= | \left\{ \textbf{n}_{i}^{j} \in N(T): \textbf{n}_{i}^{j} \neq e_i \right\}|,$$
and we need to compute this cardinality.
Taking in account the expressions for $c[i]$,
it follows that there are
\begin{itemize}
\item $\max(l(M_1),p_1)$ nodes along the first branch that are different from $e_1$;
\item$ \max(l(M_i),p_i,p_{i-1})$  nodes along the $i$-th branch  different from $e_i$, for $2 \leq i \leq \nolinebreak r-1$;
\item $\max(l(M_r),p_{r-1})$ nodes along the last branch that are different from $e_r$.
\end{itemize}
Now from  $T=T_E=(p_1,\ldots,p_{r-1})$ we deduce that the $i$-th and $i+1$-th branches have $p_i$ nodes in common for each $i=1,\ldots,r-1$.
Therefore we can conclude:
$$d(S \setminus C)= \max(l(M_1),p_1)+\sum_{k=2}^{r-1}{\max(l(M_k),p_k,p_{k-1})}+\max(l(M_r),p_{r-1})-\sum_{k=1}^{r-1}{p_k}.$$ 
Finally we have:
$$ g(S(T))=d(\mathbb{N}^r \setminus C) - d(S \setminus C)=\sum_{j=1}^{r}\left(\sum_{k=1}^{l(M_j)}{M_{j}[k]}-l(M_j)\right)+\sum_{k=1}^{r-1}{p_k},$$
and, because  $ \displaystyle \sum_{k=1}^{l(M_j)}{M_{j}[k]}-l(M_j)=g(\textrm{AS}(M_j))$, we have:
$$ g(S(T))=\sum_{k=1}^r{g(\textrm{AS}(M_k))}+\sum_{k=1}^{r-1}{p_k}, $$
and the proof is complete. \qed
\\

Now we  denote by $\textrm{Gen}(r,n)$ the set of the untwisted mulitplicity tree of all the Arf semigroups in  $\mathbb{N}^r$ with genus $n$. Given a $n \in \mathbb{N}$ we want to find a way to compute the set $\textrm{Gen}(r,n)$.
We  do that using recursion on $r$. From the previous section we know how to compute $\textrm{Gen}(1,n)$, so the base case is done.
First of all, we notice that we need $n \geq r-1$. In fact an untwisted Arf semigroup $S$ of $\mathbb{N}^r$ can be described by a tree $T=T_E=(p_1,\ldots,p_{r-1})$ with $E=\left\{M_1,\ldots,M_r\right\}$, and we have just showed that 
$$ g(S(T))=\sum_{k=1}^r{g(\textrm{AS}(M_k))}+\sum_{k=1}^{r-1}{p_k}, $$
where $g(\textrm{AS}(M_k)) \geq 0$ and $p_k \geq 1$ for all the $k$. Then $g(S(T)) \geq r-1$.

We fix a $r\geq 2$ and $n \geq r-1$ and suppose that $T=T_E=(p_1,\ldots,p_{r-1})$ is a multiplicity tree in $\textrm{Gen}(r,n)$, where $E=\left\{ M_1,\ldots,M_r \right\}$ is a collection of $r$ multiplicity sequences.

Consider $t<r$. We have:
$$ n=\sum_{j=1}^r{g(\textrm{AS}(M_j))}+\sum_{j=1}^{r-1}{p_j}=\left( \sum_{j=1}^{t}{g(\textrm{AS}(M_j))}+\sum_{j=1}^{t-1}{p_j}\right)+ p_t+\left( \sum_{j=t+1}^{r}{g(\textrm{AS}(M_j))}+\sum_{j=t+1}^{r-1}{p_j}\right), $$
therefore if we denote by 
$$ k_1=  \sum_{j=1}^{t}{g(\textrm{AS}(M_j))}+\sum_{j=1}^{t-1}{p_j} \quad \textrm{ and }  \quad   k_2=\sum_{j=t+1}^{r}{g(\textrm{AS}(M_j))}+\sum_{j=t+1}^{r-1}{p_j},$$
we have: 
$$ n-p_{t}=k_1+k_2.$$
Now, we have:
\begin{itemize} 
\item The tree $T^1=T_{E_{t}}=(p_1,\ldots,p_{t-1})$, with $E_{t}=\left\{ M_1,\ldots,M_{t} \right\}$, belongs to $\textrm{Gen}(t,k_1)$ ($k_1 \geq t-1$);
\item The tree $T^2=T_{E^{t}}=(p_{t+1},\ldots,p_{r-1})$, with $E^{t}=\left\{ M_{t+1},\ldots,M_{r} \right\}$, belongs to 
\\
$\textrm{Gen}(r-t,n-p_t-k_1)$ ;
\item $ 1 \leq p_{t} \leq \textrm{Comp}(M_{t},M_{t+1});$
\item $k_1+k_2\geq r-2$ $\Rightarrow  1 \leq p_{t} \leq n-r+2$;
\item $k_2 \geq r-1-t$  $\Rightarrow k_1 \leq n-p_t-r+1+t$.
\end{itemize}
Now, for each $ 1\leq p \leq n-r+2$ and $ t-1 \leq k(p) \leq n-p-r+1+t$ we define the set
$$ I_r^n(t,p,k(p))=\left\{ T_E=(p_1,\ldots,p_{r-1}): E=\left\{M_1,\ldots,M_r\right\}, T_{E_{t}}=(p_1,\ldots,p_{t-1}) \in \right.$$   $$ \in \textrm{Gen}(t,k(p)), \left.T_{E^{t}}=(p_{t+1},\ldots,p_{r-1}) \in \textrm{Gen}(r-t,n-p-k(p)), \right.$$ $$ \left.\textrm{ and } p_t=p\leq \textrm{Comp}(M_{t},M_{t+1}) \right\}.$$ 
So we can deduce that $T$ belongs to the following set:
$$ S(r,n)= \bigcup_{p=1}^{n-r+2}\left({ \bigcup_{k(p)=t-1}^{n-p-r+t+1}{ I_r^n(t,p,k(p))}}\right).$$
With the inverse implications, it is very easy to show that $S(r,n) \subseteq \textrm{Gen}(r,n)$.
Then we have   $ \textrm{Gen}(r,n)=S(r,n)$.

Notice that the computation of $\textrm{Gen}(r,n)$ involves the computation of $\textrm{Gen}(t,k)$ and $\textrm{Gen}(r-t,k)$, then using recursion and the base case we can solve our problem.

\begin{remark2}
If we have  $T=T_E=(p_1,\ldots,p_{r-1})$ with $E=\left\{M_1,\ldots,M_r \right\}$, we will denote by $T^{-1}=T_{E^{-1}}=(p_{r-1},\ldots,p_1)$ where $E^{-1}=\left\{ M_r,\ldots,M_1\right\}$.

It is clear that if $T \in \textrm{Gen}(r,n)$ then $T^{-1} \in \textrm{Gen}(r,n)$ too.
If $U$ is a set of multiplicity trees, we denote by $U^{-1}=\left\{ T^{-1}: T \in U \right\}$. 
\end{remark2}We have the following proposition:

\begin{raff} \label{raff}

If $t \leq \frac{r}{2}$ and $\lambda(p) \geq \frac{n-p-1}{2}$ then
$$ \textrm{Gen}(r,n)= \bigcup_{p=1}^{n-r+2}\left({ \bigcup_{k(p)=t-1}^{\lambda(p)}\left({ I_r^n(t,p,k(p)) \cup \left(I_r^n(t,p,k(p))\right)^{-1}}\right)} \right).$$
\end{raff}
\noindent \textbf{Proof}.
Consider $t \leq \frac{r}{2}$ and $\lambda(p) \geq \frac{n-p-1}{2}$, we show that for any $p=1,\ldots,n-r+2$, we have
$$ \bigcup_{k(p)=\lambda(p)+1}^{n-p-r+t+1}{ I_r^n(t,p,k(p))} \subseteq \left(\bigcup_{k(p)=t-1}^{\lambda(p)}{ I_r^n(t,p,k(p))}\right)^{-1}.$$
Consider $\displaystyle T \in \bigcup_{k(p)=\lambda(p)+1}^{n-p-r+t+1}{ I_r^n(t,p,k(p))}$, then 
$$ T=T_E=(p_1,\ldots,p_{r-1}) \textrm{ with } E=\left\{ M_1,\ldots,M_r\right\} \textrm{ and } T_{E^t}=(p_{t+1},\ldots,p_{r-1}) \in \textrm{Gen}(r-t,n-p-\overline{k}),$$
where $\lambda(p)+1 \leq \overline{k} \leq n-p-r+t+1$.

We want to show that $$ \displaystyle T \in  \left(\bigcup_{k(p)=t-1}^{\lambda(p)}{ I_r^n(t,p,k(p))}\right)^{-1}  \iff T^{-1} \in \bigcup_{k(p)=t-1}^{\lambda(p)}{ I_r^n(t,p,k(p))}.$$
Let us consider the subtree $T^1$ of $T^{-1}$ with  $T^1=T_{(E^{-1})_t}=(p_{r-1},\ldots,p_{r-t+1})$, where $(E^{-1})_t=\left\{M_r,\ldots,M_{r-t+1} \right\}$.

Now we have:
$$ t \leq \frac{r}{2} \Rightarrow 2t \leq r \Rightarrow t\leq r-t \Rightarrow t+1 \leq r-t+1,$$
and from this and from $ T_{E^t}=(p_{t+1},\ldots,p_{r-1}) \in \textrm{Gen}(r-t,n-p-\overline{k})$ we can deduce that
$$ T^1=T_{(E^{-1})_t}=(p_{r-1},\ldots,p_{r-t+1}) \in \textrm{Gen}(t,x),$$
where $t-1 \leq x \leq n-p-\overline{k}$.

Thus, from  $\lambda(p)+1 \leq \overline{k}$ and from $\lambda(p) \geq \frac{n-p-1}{2}$ we deduce:
$$  -\overline{k} \leq -\lambda(p)-1 \textrm{ and } n-p-1 \leq 2\lambda(p) \Rightarrow x \leq n-p-\lambda(p)-1 \leq \lambda(p).$$ 
Therefore, from $T^{-1} \in \textrm{Gen}(r,n)$ and from the previous inequality we  have 
$$ T^{-1} \in  I_r^n(t,p,x) \subseteq \bigcup_{k(p)=t-1}^{\lambda(p)}{ I_r^n(t,p,k(p))},$$
and the claim is proved.
Therefore we have:
$$ \textrm{Gen}(r,n) \supseteq \bigcup_{p=1}^{n-r+2}\left({ \bigcup_{k(p)=t-1}^{\lambda(p)}\left({ I_r^n(t,p,k(p)) \cup \left(I_r^n(t,p,k(p))\right)^{-1}}\right)} \right) \supseteq $$ $$ \supseteq  \bigcup_{p=1}^{n-r+2}\left({ \bigcup_{k(p)=t-1}^{n-p-r+t+1}{ I_r^n(t,p,k(p))}}\right)=\textrm{Gen}(r,n),$$
and the proof is complete. \qed

The previous proposition suggest us an easier way to compute $\textrm{Gen}(r,n)$. In fact we have to consider a smaller amount of sets of the type $I_r^n(t,p,k(p))$, completing the computation with sets of the type $\left(I_r^n(t,p,k(p))\right)^{-1}$ which are very easy to obtain once we have $I_r^n(t,p,k(p))$.
To speed up the process is also useful to set $\lambda(p)=\left \lceil{\frac{n-p-1}{2}}\right \rceil $ and $t=\left \lfloor{\frac{r}{2}}\right \rfloor$.

Denote by $\overline{\textrm{Gen}(r,n)}$ the set of all the multiplicity trees (twisted and untwisted) of the Arf semigroups in $\mathbb{N}^r$ with genus $n$. We are already able to compute all the untwisted ones. Suppose  that $T$ is a twisted tree in 
 $\overline{\textrm{Gen}(r,n)}$. Then there exists a permutation $\sigma \in S^r$ such that $\sigma(T)$ is untwisted. From the formula of the genus it is very easy to see that $\sigma(T) \in \textrm{Gen}(r,n)$. Thus we have:
$$ \overline{\textrm{Gen}(r,n)}=\bigcup_{\sigma \in S^r} \left\{ \sigma^{-1}(T): T \in \textrm{Gen}(r,n)  \right\}.$$
\begin{example6}
We compute $\textrm{Gen}(2,3)$. In this case we have $t=1$.

We have to consider the sets $I_2^3(1,p,k(p))$ and $\left(I_2^3(1,p,k(p))\right)^{-1}$  for  each $ 1\leq p \leq 3$ and $ 0 \leq k(p) \leq \left \lceil{\frac{3-p-1}{2}}\right \rceil $.

\textbf{Case:} $p=1$  $\Rightarrow 0\leq k(1) \leq 1$
\begin{itemize}
\item $k(1)=0$.
To compute $I_2^3(1,1,0)$ we need $M_1 \in \textrm{Gen}(1,0)$ and $M_2 \in \textrm{Gen}(1,2)$. They are 
$$  \textrm{Gen}(1,0)=\left\{ [1] \right\}  \textrm{ and }  \textrm{Gen}(1,2)=\left\{ [3],[2,2] \right\}.$$
Thus we can consider  \begin{itemize} \item $E_1=\left\{M_1=[1],M_2=[3]\right\}$. 

In this case $1=p \leq \textrm{Comp}(M_1,M_2)=2$ then the tree $T_1=T_{E_1}=(1)$ belongs to   $\textrm{Gen}(2,3).$
\item $E_2=\left\{M_1=[1],M_2=[2,2]\right\}$. 

In this case $1=p \leq \textrm{Comp}(M_1,M_2)=3$ then the tree $T_2=T_{E_2}=(1)$ belongs to   $\textrm{Gen}(2,3).$
 \end{itemize}
Therefore $I_2^3(1,1,0)=\left\{ T_1,T_2\right\}.$
We have now to compute $\left(I_2^3(1,1,0)\right)^{-1}$
\begin{itemize}
\item $T_3=T_1^{-1}=T_{E_3}=(1) \in \textrm{Gen}(2,3),$ with $E_3=E_1^{-1}=\left\{ [3],[1]\right\}.$
\item $T_4=T_2^{-1}=T_{E_4}=(1) \in \textrm{Gen}(2,3),$ with $E_4=E_2^{-1}=\left\{ [2,2],[1]\right\}.$
\end{itemize}
Therefore $\left(I_2^3(1,1,0)\right)^{-1}=\left\{ T_3,T_4\right\}.$
\item $k(1)=1$.
To compute $I_2^3(1,1,1)$ we need $M_1 \in \textrm{Gen}(1,1)$ and $M_2 \in \textrm{Gen}(1,1)$.  We have 
$$  \textrm{Gen}(1,1)=\left\{ [2] \right\}. $$
Thus we only have to  consider  \begin{itemize} \item $E_5=\left\{M_1=[2],M_2=[2]\right\}$. 

In this case $1=p \leq \textrm{Comp}(M_1,M_2)=+\infty$ then the tree $T_5=T_{E_5}=(1)$ belongs to   $\textrm{Gen}(2,3).$
 \end{itemize}
Therefore $I_2^3(1,1,1)=\left\{ T_5 \right\}$. In this case $I_2^3(1,1,1)=\left(I_2^3(1,1,1) \right)^{-1}$.
\end{itemize}
\textbf{Case:} $p=2$  $\Rightarrow k(2)=0$
\begin{itemize}
\item $k(2)=0$.
To compute $I_2^3(1,2,0)$ we need $M_1 \in \textrm{Gen}(1,0)$ and $M_2 \in \textrm{Gen}(1,1)$. 

Thus we only have to consider  \begin{itemize} \item $E_6=\left\{M_1=[1],M_2=[2]\right\}$. 

In this case $2=p \leq \textrm{Comp}(M_1,M_2)=2$ then the tree $T_6=T_{E_6}=(2)$ belongs to   $\textrm{Gen}(2,3).$
 \end{itemize}
Therefore $I_2^3(1,2,0)=\left\{T_6\right\}.$
We have now to compute $\left(I_2^3(1,2,0)\right)^{-1}$
\begin{itemize}
\item $T_7=T_6^{-1}=T_{E_7}=(2) \in \textrm{Gen}(2,3),$ with $E_7=E_6^{-1}=\left\{ [2],[1]\right\}.$
\end{itemize}
Therefore $\left(I_2^3(1,2,0)\right)^{-1}=\left\{ T_7\right\}.$
 \end{itemize}
\textbf{Case:} $p=3$  $\Rightarrow k(3)=0$
\begin{itemize}
\item $k(3)=0$.

To compute $I_2^3(1,3,0)$ we need $M_1 \in \textrm{Gen}(1,0)$ and $M_2 \in \textrm{Gen}(1,0)$. 

Thus we only have to consider   \begin{itemize} \item $E_8=\left\{M_1=[1],M_2=[1]\right\}$. 

In this case $3=p \leq \textrm{Comp}(M_1,M_2)=+\infty$ then the tree $T_8=T_{E_8}=(3)$ belongs to   $\textrm{Gen}(2,3).$ \end{itemize} 
Therefore $I_2^3(1,3,0)=\left\{ T_8 \right\}$. In this case $I_2^3(1,3,0)=\left(I_2^3(1,3,0) \right)^{-1}$.
Thus $\textrm{Gen}(2,3)=\left\{ T_1, T_2, T_3, T_4, T_5, T_6, T_7, T_8\right\}$. We obviously have $\textrm{Gen}(2,3)=\overline{\textrm{Gen}(2,3)}$ because we have only two branches.
\end{itemize}
\vskip 0.3in

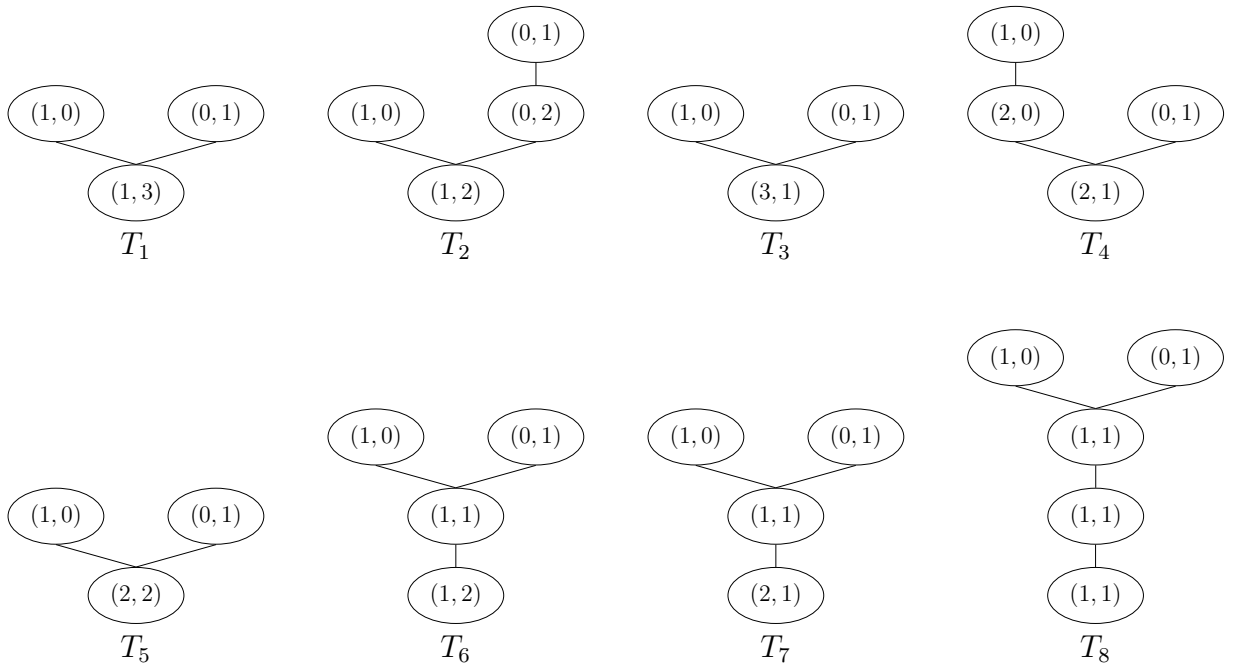
\begin{tikzpicture}[grow'=up,sibling distance=32pt,scale=.75]
\tikzset{level distance=40pt,every tree node/.style={draw,ellipse}} \Tree [ .$(1,3)$  [ .$(1,0)$  ] [ .$(0,1)$ ]  ] ]\node[below]at(current bounding box.south){$T_1$}; \end{tikzpicture}
\hskip 0.3in
\begin{tikzpicture}[grow'=up,sibling distance=32pt,scale=.75]
\tikzset{level distance=40pt,every tree node/.style={draw,ellipse}} \Tree [ .$(1,2)$  [ .$(1,0)$  ] [ .$(0,2)$  [.$(0,1)$ ]]  ]\node[below]at(current bounding box.south){$T_2$}; \end{tikzpicture}
\hskip 0.3in
\begin{tikzpicture}[grow'=up,sibling distance=32pt,scale=.75]
\tikzset{level distance=40pt,every tree node/.style={draw,ellipse}} \Tree [ .$(3,1)$  [ .$(1,0)$  ] [ .$(0,1)$ ]  ] ]\node[below]at(current bounding box.south){$T_3$}; \end{tikzpicture}
\hskip 0.3in
\begin{tikzpicture}[grow'=up,sibling distance=32pt,scale=.75]
\tikzset{level distance=40pt,every tree node/.style={draw,ellipse}} \Tree [ .$(2,1)$  [ .$(2,0)$ [.$(1,0)$ ] ] [ .$(0,1)$ ]  ]\node[below]at(current bounding box.south){$T_4$}; \end{tikzpicture}
\vskip 0.3in
\begin{tikzpicture}[grow'=up,sibling distance=32pt,scale=.75]
\tikzset{level distance=40pt,every tree node/.style={draw,ellipse}} \Tree [ .$(2,2)$  [ .$(1,0)$  ] [ .$(0,1)$ ]  ] ]\node[below]at(current bounding box.south){$T_5$}; \end{tikzpicture}
\hskip 0.3in
\begin{tikzpicture}[grow'=up,sibling distance=32pt,scale=.75]
\tikzset{level distance=40pt,every tree node/.style={draw,ellipse}} \Tree [ .$(1,2)$ [ .$(1,1)$  [ .$(1,0)$  ] [ .$(0,1)$ ]  ] ] ]\node[below]at(current bounding box.south){$T_6$}; \end{tikzpicture}
\hskip 0.3in
\begin{tikzpicture}[grow'=up,sibling distance=32pt,scale=.75]
\tikzset{level distance=40pt,every tree node/.style={draw,ellipse}} \Tree [ .$(2,1)$ [ .$(1,1)$  [ .$(1,0)$  ] [ .$(0,1)$ ]  ] ] ]\node[below]at(current bounding box.south){$T_7$}; \end{tikzpicture}
\hskip 0.3in
\begin{tikzpicture}[grow'=up,sibling distance=32pt,scale=.75]
\tikzset{level distance=40pt,every tree node/.style={draw,ellipse}} \Tree [ .$(1,1)$ [ .$(1,1)$ [ .$(1,1)$  [ .$(1,0)$  ] [ .$(0,1)$ ]  ] ] ] ]\node[below]at(current bounding box.south){$T_8$}; \end{tikzpicture}

\end{example6}
In the following table we report the cardinality of $\textrm{Gen}(2,n)$ for $n$ up to $32$.
\begin{center}
	\resizebox{16cm}{3cm}{
		\begin{tabular}{ || c | c || c |c|| c| c|| c| c||}
			\hline \hline
			$n$ & $| \textrm{Gen}(2,n)|$& $n$ & $|  \textrm{Gen}(2,n)|$& $n$ & $|  \textrm{Gen}(2,n)|$& $n$ & $| \textrm{Gen}(2,n)|$ \\ 
			
			\hline \hline
			$1$ & 1 & $9$ & 251 & $17$ & 4386& $25$ & 35203 \\
			\hline
			$2$ & 3 & $10$ & 385 & $18$ & 5874& $26$ & 44209 \\
			\hline
			$3$ &8 & $11$ & 577 & $19$ &7773& $27$ & 55175 \\
			\hline
			$4$ & 16 & $12$ & 837 & $20$ & 10195& $28$ & 68493 \\
			\hline
			$5$ & 32 & $13$ & 1207 & $21$ & 13270& $29$ & 84540 \\
			\hline
			$6$ & 56 & $14$ & 1701 & $22$ & 17138& $30$ & 103898 \\
			\hline
			$7$ & 99 & $15$ & 2361 & $23$ & 21922& $31$ & 127031 \\
			\hline
			$8$ & 157 & $16$ & 3239 & $24$ & 27882& $32$ & 154681 \\
			\hline \hline
	\end{tabular}}
\end{center}

Using the previous results, it is easy to implement an algorithm that computes the number of  untwisted Arf semigroups of $ \mathbb{N}^r$ with a given genus $n$. Each entry of the following table is such a number, for $ 1 \leq r \leq 16$ and $0\leq n \leq 15$.
\begin{center}
	\resizebox{16cm}{!}{
		\begin{tabular}{| | c || c| c | c | c| c| c| c| c| c| c| c| c| c| c| c|c ||p{5cm} }
			\hline \hline 
			$r \setminus n $ & \textbf{0}& \textbf{1} &\textbf{2}  & \textbf{3} &\textbf{4} & \textbf{5}&\textbf{6}&\textbf{7}&\textbf{8}&\textbf{9}&\textbf{10}&\textbf{11}&\textbf{12}&\textbf{13}&\textbf{14}&\textbf{15}\\ \hline \hline
			\textbf{1} & 1&1 &2  & 3 &4 & 6&8&10&13&17&21&26&31&36&47&55  \\ \hline
			\textbf{2}& 0&1 &3  & 8 &16 & 32&56&99&157&251&385&577&837&1207&1701&2361  \\  \hline
			\textbf{3} & 0&0 &1  & 5 &18 & 49&120&263&543&1048&1943&3458&5957&9957&16246&25896 \\
			\hline
			\textbf{4} & 0&0 &0  & 1 &7 & 32&110&324&846&2032&4544&9620&19420&37686&70618&128399 \\
			\hline
			\textbf{5} & 0&0 &0  &0 &1 & 9&50&207&716&2169&5958&15119&35994&81196&175001&362501  \\
			\hline
			\textbf{6} & 0&0 &0  & 0 &0 & 1&11&72&348&1384&4772&14769&41919&110859&276257&654422   \\
			\hline
			\textbf{7}  & 0&0 &0  & 0 &0 & 0&1&13&98&541&2432&9403&32385&101658&295681&806530  \\
			\hline
			\textbf{8}  & 0&0 &0  & 0 &0 & 0&0&1&15&128&794&3980&17050&64678&222474&705806\\
			\hline
			
			\textbf{9}  & 0&0 &0  & 0 &0 & 0&0&0&1&17&162&1115&6164&28973&120016&448873 \\
			\hline
			\textbf{10}  & 0&0 &0  & 0 &0 & 0&0&0&0&1&19&200&1512&9136&46736&209871\\
			\hline
			\textbf{11}  & 0&0 &0 & 0 &0 & 0&0&0&0&0&1&21&242&1993&13064&72239\\
			\hline
			\textbf{12}  & 0&0 &0  & 0 &0 & 0&0&0&0&0&0&1&23&288&2566&18132 \\
			\hline
			\textbf{13} & 0&0 &0  & 0 &0 & 0&0&0&0&0&0&0&1&25&338&3239 \\
			\hline
			\textbf{14}  & 0&0 &0  & 0 &0 & 0&0&0&0&0&0&0&0&1&27&392\\
			\hline
			\textbf{15}  & 0&0 &0  & 0 &0 & 0&0&0&0&0&0&0&0&0&1&29\\
			\hline
			\textbf{16} & 0&0 &0  & 0 &0 & 0&0&0&0&0&0&0&0&0&0&1 \\
			
			\hline \hline
			\end{tabular}}
		\end{center}
From the previous table, we can also deduce, considering the sum of all the entries in the corresponding column, the number of all the local untwisted Arf semigroups with a given genus $n$ (in all the possible dimensions). We call $$\textrm{NG}(n)=| \left\{ S \textrm{Arf semigroup}: g(S)=n \right\}|, $$ such a number.

Thus we have
\begin{center}
	\resizebox{17cm}{!}{
		\begin{tabular}{| | c || c| c | c | c| c| c| c| c| c| c| c| c| c| c| c|c ||p{5cm} }
			\hline \hline 
			$n $ & \textbf{0}& \textbf{1} &\textbf{2}  & \textbf{3} &\textbf{4} & \textbf{5}&\textbf{6}&\textbf{7}&\textbf{8}&\textbf{9}&\textbf{10}&\textbf{11}&\textbf{12}&\textbf{13}&\textbf{14}&\textbf{15}\\ \hline 
			$\textrm{NG}(n)$ & 1&2 &6  & 17 &46 & 129&356&989&2737&7588&21031&58289&161535&447693&1240773&3438746  \\ \hline
			 \hline
	\end{tabular}}
\end{center}

\begin{example8}
In the following table we report the cardinality of the sets $\overline{\textrm{Gen}(r,n)}$ for $ 1 \leq r \leq 9$ and $0\leq n \leq 8$.
\begin{center}
	\resizebox{!}{4cm}{
		\begin{tabular}{| | c || c| c | c | c| c| c| c| c| c|| p{5cm} }
			\hline \hline 
			$r \setminus n $ & \textbf{0}& \textbf{1} &\textbf{2}  & \textbf{3} &\textbf{4} & \textbf{5}&\textbf{6}&\textbf{7}&\textbf{8}\\ \hline \hline
			\textbf{1} & 1&1 &2  & 3 &4 & 6&8&10&13  \\ \hline
			\textbf{2} & 0&1 &3  & 8 &16 & 32&56&99&157   \\  \hline
			\textbf{3} & 0&0 &1  & 6 &22 & 61&151&334&693  \\
			\hline
			\textbf{4}  & 0&0 &0  & 1 &10 &51&189&576&1555  \\
			\hline
			\textbf{5}  & 0&0 &0  & 0 &1 & 15&105&505&1906  \\
			\hline
			\textbf{6}  & 0&0 &0  & 0 &0 & 1&21&197&1208    \\
			\hline
			\textbf{7}  & 0&0 &0  & 0 &0 & 0&1&28&343   \\
			\hline
			\textbf{8}   & 0&0 &0  & 0 &0 & 0&0&1&36 \\
			\hline
			
			\textbf{9}   & 0&0 &0  &0 &0 &0&0&0&1  \\
			\hline
			 \hline
	\end{tabular}}
\end{center}
\end{example8}
\begin{otherlanguage}{english}

\end{otherlanguage}

E-mail address: giuseppezito@hotmail.it

\end{document}